\documentclass[11pt]{article}
\usepackage{amssymb}

\advance \textheight by -1 \topmargin
\advance \textheight by -1 \headheight
\advance \textheight by -1 \headsep
\advance \textheight by -2 \footskip
\vsize = \textheight
\hsize = \textwidth
\newcommand{\BFW}{\mbox{\rm BFW}}
\newcommand{\Per}{\mbox{\rm Per}}
\newcommand{\Mat}{\mbox{\rm Mat}}
\newcommand{\Sym}{\mbox{\rm Sym}}
\newcommand{\diag}{\mbox{\rm diag}}
\newcommand{\GL}{\mbox{\rm GL}}
\newcommand{\Stab}{\mbox{\rm Stab}}
\newcommand{\Aut}{\mbox{\rm Aut}}
\newcommand{\End}{\mbox{\rm End}}
\newcommand{\trace}{\mbox{\rm trace}}
\newcommand{\Sp}{\mbox{\rm Sp}}
\newtheorem{theorem}{Theorem}[section]

\newtheorem{Conjecture}[theorem]{Conjecture}

\newtheorem{remark}[theorem]{Remark}

\newtheorem{kor}[theorem]{Corollary}
\newtheorem{defi}[theorem]{Definition}

\newtheorem{lemma}[theorem]{Lemma}

\renewcommand{\setminus}{-}
\newcommand{\fa}{\mbox{ for}\mbox{ all }}
\newcommand{\be}{\begin{enumerate}}
\newcommand{\ee}{\end{enumerate}}
\newcommand{\bi}{\begin{itemize}}
\newcommand{\ei}{\end{itemize}}

\newcommand{\ba}{\begin{array}}
\newcommand{\ea}{\end{array}}

\newcommand{\ra}{\rightarrow}

\newcommand{\teiltnicht}{\mbox{ $\! \mid \!\!\! \not \ $ }}

\newcommand{\Z}{{\mathbb{Z}}}
\newcommand{\Q}{{\mathbb{Q}}}

\newcommand{\N}{{\mathbb{N}}}
\newcommand{\R}{{\mathbb{R}}}
\newcommand{\A}{{\mathbb{A}}}

\newcommand{\C}{{\mathbb{C}}}

\renewcommand{\em}{\sf}

\begin{document}

\Large
\begin{center}
{\bf Hecke actions on certain strongly modular genera of lattices}
\end{center}
\normalsize
\begin{center}
{\it Gabriele Nebe} 
 and {\it Maria Teider} 
\end{center}

\small 
{\sc Abstract}
We calculate the action of some 
Hecke operators on spaces of modular forms spanned
by the Siegel theta-series of certain genera of strongly 
modular lattices closely related to the Leech lattice.
Their eigenforms provide explicit examples of Siegel cusp forms.

\normalsize

\section{Introduction}

One of the most remarkable lattices in Euclidean space is the
Leech lattice, the unique even unimodular lattice
$\Gamma _1\subset (\R ^{24}, (,)) $ of dimension 24 that
does not contain vectors of square length 2.
Here a lattice $\Lambda \subset (\R^n,(,))$ is called 
{\em unimodular}, if $\Lambda $ equals its {\em dual lattice}
$$\Lambda ^{\#} := \{ x\in \R ^n \mid (x,\lambda ) \in \Z \fa \lambda \in \Lambda \} $$
and {\em even}, if the quadratic form $x\mapsto (x,x) $ takes only even values 
on $\Lambda $.
\cite{SieMod} studies spaces of Siegel modular forms generated by the Siegel
theta-series of the 24 isometry classes of lattices in the genus of
$\Gamma _1$. 
The present paper extends this investigation to further genera of lattices,
closely related to $\Gamma _1$.
A unified construction is given in \cite{RS98s}:
Consider the Matthieu group 
$M_{23} \leq \Aut (\Gamma _1)$, where the {\em automorphism group} of 
a lattice $\Lambda \subset (\R ^n,(,)) $ is
$\Aut(\Lambda ):= \{ g\in O(n)  \mid  \Lambda g = \Lambda \} .$
Let $g\in M_{23}$ be an element of square-free order 
$l:=|\langle g \rangle |$. Then 
$$l\in \{ 1,2,3,5,6,7,11,14,15,23 \} =: {\cal N} = \{ n\in \N \mid 
\sigma_1(n):=\sum _{d \mid n} d \mbox{ divides } 24 \} $$
and for each $l\in {\cal N}$, there is an up to conjugacy  unique
cyclic subgroup $\langle g \rangle \leq M_{23}$ of order $l$.
Let $\Gamma _l := \{ \lambda \in \Gamma _1 \mid \lambda g = \lambda \}$
denote the fixed lattice of $g$.
Then $\Gamma _l$ is an extremal strongly modular lattice of level $l$ and
of dimension $2k_l$, where $$k_l:=12 \sigma_0(l)/\sigma _1(l)$$ 
and $\sigma_0(l)$ denotes the number of divisors of $l$.
In particular $\Gamma _1$ is the Leech lattice,
$\Gamma _2$ the 16-dimensional Barnes-Wall lattice and
$\Gamma _3$ the Coxeter-Todd lattice of dimension 12.

Let $\Lambda $ be an even lattice.
The minimal $l\in \N $ for which $\sqrt{l} \Lambda ^{\# }$ is even,
is called the {\em level} of $\Lambda $.
Then $l\Lambda ^{\#} \subset \Lambda $.
For an exact divisor $d$ of $l$ let 
$$\Lambda ^{\#,d} := \Lambda ^{\#} \cap \frac{1}{d} \Lambda $$
denote the {\em $d$-partial dual} of $\Lambda $.
A lattice $\Lambda $ is called {\em strongly $l$-modular}, if
$\Lambda $ is isometric to $\sqrt{d} \Lambda ^{\# ,d}$ for all 
exact divisors $d$ of the level $l$ of $\Lambda $.
If $l$ is a prime, this coincides with the notion of {\em modular}
lattices, which just means that the lattice is similar to its dual
lattice.
The Siegel theta-series 
$$\Theta ^{(m)} _{\Lambda } (Z) := \sum _{(\lambda _1,\ldots , \lambda _m) \in 
\Lambda ^m }\exp ( i \pi  \trace ((\lambda _i,\lambda _j)_{i,j} Z)) $$
(which is a holomorphic function on
the Siegel halfspace ${\cal H}^{(m)} = $ $ \{ Z\in \Sym _m(\C)
\mid $ $ \Im (Z) $ positive  definite $ \} $) 
of a strongly $l$-modular lattice is a modular form for the
$l$-th congruence subgroup
$\Gamma _0^{(m)}(l)$ of $ \Sp_{2m}(\Z )$ (to a certain character)
invariant under 
all Atkin-Lehner-involutions (cf. \cite{Andrianov}).
In particular for $m=1$ and $l\in {\cal N}$ 
the relevant ring of modular forms is a polynomial
ring in 2 generators as shown in \cite{Quebmodular}, \cite{Quebstmodular}.
Explicit generators of this ring allow to bound the 
minimum of an $n$-dimensional strongly 
$l$-modular lattice $\Lambda $ with $l\in {\cal N}$,
$$
\min (\Lambda ) := \min _{0\neq \lambda \in \Lambda }(\lambda, \lambda )
\leq 2 + 2 \lfloor  \frac{n}{2k_l} \rfloor .$$
Lattices $\Lambda $ achieving this bound 
are called {\em extremal}.
For all $l \in {\cal N}$ there is a unique 
extremal strongly $l$-modular lattices of dimension $2k_l$
and this is
the lattice $\Gamma _l$ described above.
All the genera are presented in the nice survey article
\cite{SchaSchuPi}.

In this paper we investigate the spaces of Siegel modular forms
generated by the Siegel theta-series of the lattices 
in the genus ${\cal G} (\Gamma _l)$ 
for $l\in {\cal N}$ using similar methods as for
the case $l=1$ which is treated in \cite{SieMod}.
The vector space ${\cal V} := {\cal V}({\cal G})$ 
of all complex formal linear combinations of the
isometry classes of lattices in 
any genus ${\cal G}$  forms a  finite dimensional
commutative
$\C $-algebra  with positive definite Hermitian scalar product.
Taking theta-series defines linear operators 
$\Theta ^{(m)}$ from ${\cal V}$ into a certain space of
modular forms and hence a filtration of ${\cal V}$ by
the kernels of these operators. 
This filtration behaves nicely under the multiplication 
and is invariant under all Hecke-operators.
With the Kneser neighbouring process we construct a 
family of commuting self-adjoint linear operators on ${\cal V}$.
Their common eigenvectors provide
explicit examples of Siegel cusp forms.

The genera ${\cal G} (\Gamma _l) $ ($l\in {\cal N}$) 
share the following properties:

\begin{kor}
Let $l\in {\cal N}$ and
let $p$ be the smallest prime not dividing  $l$.
The mapping $\Theta ^{(k_l)}$ is injective on ${\cal V}({\cal G}(\Gamma _l))$.
For $l\neq 7$, the construction described in \cite{BFW} (see Paragraph \ref{BFW})
gives a non-zero cusp form $\BFW (\Gamma _l,p) =
\Theta ^{(k_l)}(\Per (\Gamma _l,p))$.
The eigenvalue of the Kneser operator $K_2$ at
the eigenvector $\Per (\Gamma _l,p) $
is the negative of the number of pairs of
minimal vectors in $\Gamma _l$ which is also the
minimal eigenvalue of $K_2$.
\end{kor}

\begin{remark}
In Section \ref{results} we also list the eigenvalues of some of the
operators $T(q)$ defined in Subsection \ref{Hecke}.
These eigenvalues suggest that for even values of $k_l$, the
cusp form $\BFW (\Gamma _l,p)$ is a generalized Duke-Imamoglu-Ikeda lift 
(see \cite{Ikeda}) of the elliptic cusp form of minimal weight
$k_l$. 
\end{remark}

{\bf Acknowledgement.} 
We thank R. Schulze-Pillot for helpful comments, suggestions and 
references.

\section{Methods}

The general method has already been explained in \cite{SieMod}
(see also \cite{SchuPi0},\cite{SchuPi1}, 
\cite{SchuPi2} and \cite{Birch} for similar strategies).

\subsection{The algebra ${\cal V} = {\cal V}({\cal G})$}

Let ${\cal G}$ be a genus of lattices in the Euclidean space 
$(\R ^{2k},(,))$.
Then ${\cal G}$ is the disjoint union of finitely many isometry classes
$${\cal G} =  [\Lambda _1] \cup \ldots \cup [\Lambda _h ] .$$

Let ${\cal V} := {\cal V}({\cal G})\cong \C ^h$ be the complex vector space  with basis
($[\Lambda _1], \ldots , [\Lambda _h])$.
Let ${\cal V}_{\Q } = \langle 
[\Lambda _1], \ldots , [\Lambda _h]  \rangle _{\Q } \cong \Q ^h$ 
be the rational span of the basis.

The space ${\cal V}$ can be identified with the algebra ${\cal A} $ 
of complex
functions on the double cosets 
$G(\Q) \backslash G(\A) / \Stab _{G(\A )} (\Lambda _{\A } ) 
= \cup _{i=1}^h G(\Q) x_i \Stab _{G(\A )} (\Lambda  _{\A }) $
where $G$ is the integral form of the real orthogonal group $G(\R )= O_{2k}$
defined by $\Lambda _1$, $\A $ denotes the ring
of rational ad\`eles and $\Lambda _{\A} $ the ad\'elic completion of $\Lambda _1$. 
If $\chi _{i}$ denotes the characteristic function 
mapping 
$G(\Q) x_j \Stab _{G(\A )} (\Lambda  _{\A }) $ to $\delta _{ij}$  and 
$\Lambda _i = x_i \Lambda _1$  ($i=1,\ldots , h$) then 
the isomorphism maps $[\Lambda _i ]$ 
to $|\Aut (\Lambda _i)| \chi _i $.
The usual Petersson scalar product then translates into the 
 Hermitian scalar product on ${\cal V}$ defined by
$$< [\Lambda _i], [\Lambda _j] > := \delta _{ij} |\Aut (\Lambda _i)| $$
and the multiplication of ${\cal A}$ defines 
a commutative and associative multiplication $\circ $ on ${\cal V}$ with
$$[ \Lambda _i ] \circ [ \Lambda_j ] := \# (\Aut(\Lambda _i )) \delta_{i,j}
[ \Lambda _i ]$$
(see for instance  \cite[Section 1.1]{Boe}).
Note that the Hermitian form
$\langle , \rangle $ is associative, i.e.
$$\langle v_1 \circ v_2 , v_3 \rangle = \langle v_1 , v_2 \circ v_3 \rangle
\fa v_1,v_2,v_3\in {\cal V} . $$

\subsection{The two basic filtrations of ${\cal V}$}

For simplicity we now assume that ${\cal G}$ consists of 
even lattices.
Let $l$ be the level of the lattices in ${\cal G}$.
Taking the degree-$n$ Siegel theta-series
$\Theta _{\Lambda _i }^{(n)} $ ($n=0,1,2,\ldots $)
of the lattices $\Lambda _i$
($i=1,\ldots , h $) then defines a linear map
$$ 
 \Theta ^{(n)} : {\cal V}\ra  M_{n,k}(l) \mbox{ by }
 \Theta ^{(n)}(\sum _{i=1 }^h
 c_i [ \Lambda _i] ):=
 \sum _{i=1}^h c_{i } \Theta _{\Lambda _i }^{(n)} $$
 with values in a space of modular forms 
 of degree $n$ and weight $k$ for the group $\Gamma _{0}^{(n)}(l)$
 (see \cite{Andrianov}).

For $n=0,\ldots, 2k $ let 
${\cal V}_{n}:= \ker (\Theta ^{(n)})$ be the kernel of this linear map.
Then we get the filtration
$${\cal V}=:{\cal V}_{-1} \supseteq {\cal V}_0 \supseteq {\cal V}_1 \supseteq \ldots \supseteq {\cal V}_{2k} = \{ 0 \}$$ where
${\cal V}_0 = \{ v =\sum _{i=1}^h c_{i } [\Lambda _i ] \mid
\sum _{i=1}^h c_{i } = 0 \}$ is of codimension 1 in ${\cal V}$.

Clearly $\Theta ^{(n)}({\cal V}_{n-1})$ is the kernel of the
Siegel $\Phi $-operator mapping $\Theta ^{(n)}({\cal V})$
onto $\Theta ^{(n-1)} ({\cal V})$.
For square-free level one even has

\begin{theorem} (see \cite[Theorem 8.1]{BoeSchuPi}) {\label{cusp}}
If $l$ is square-free, then 
$\Theta ^{(n)} ({\cal V}_{n-1}) $ is the space of cusp forms in
$\Theta ^{(n)} ({\cal V})$.
\end{theorem}

Let ${\cal W}_n := {\cal V}_n^{\perp }$ be the orthogonal complement of ${\cal V}_n$.
We then have the ascending filtration
$$0={\cal W}_{-1} \subseteq {\cal W}_0 \subseteq {\cal W}_1 \subseteq \ldots \subseteq {\cal W}_{2k} = {\cal V}.$$

By \cite[Proposition 2.3, Corollary 2.4]{SieMod}  
one has the following lemma:

\begin{lemma}\label{mult}
$${\cal W}_n \circ {\cal W}_m \subset {\cal W}_{n+m} \fa m,n \in \{ -1,\ldots , 2k \} $$
and 
 $${\cal W}_n \circ {\cal V}_m \subset {\cal V}_{m-n} \fa m>n \in \{ -1,\ldots , 2k \} .$$
 \end{lemma}

Since theta-series have rational coefficients, 
both filtrations are rational, i.e. 
${\cal V}_n = \C \otimes ({\cal V}_n \cap {\cal V}_{\Q })$ 
and ${\cal W}_n = \C \otimes ({\cal W}_n \cap {\cal V}_{\Q })$, hence the same statements
hold when ${\cal V}$ is replaced by ${\cal V}_{\Q }$.

\subsection{The Borcherds-Freitag-Weissauer cusp form}{\label{BFW}}

The article \cite{BFW} gives a quite general construction of a
cusp form of degree $k$. 
Let $\Lambda  $ be a $2k$-dimensional even lattice and
choose some prime $p$ such that the quadratic space
$(\Lambda / p\Lambda , Q_p) $  (where $Q_p(x) := \frac{1}{2} (x,x) + p\Z$)
is isometric to
the sum of $k$ hyperbolic planes.
Fix a totally isotropic subspace 
$F$ of $\Lambda / p\Lambda $ of dimension $k$.
For $\lambda := (\lambda _1,\ldots, \lambda _k)  \in \Lambda ^k$
we put 
$E(\lambda ):=\langle \lambda _1,\ldots, \lambda _k \rangle + p\Lambda $
and 
$S(\lambda ) := \frac{1}{p} ((\lambda _i, \lambda _j )_{i,j}) \in 
\Sym _k (\R )$.
Define $ \epsilon (E(\lambda) ) = \epsilon (\lambda ):= (-1) ^{\dim (F\cap E(\lambda ))} $
if $E(\lambda )$ is a $k$-dimensional totally isotropic subspace of 
$\Lambda /p \Lambda $ and 
$\epsilon (E(\lambda) ) = \epsilon (\lambda ) := 0$  otherwise.

\begin{defi}
$\BFW(\Lambda ,p) (Z)  := \sum _{\lambda \in \Lambda ^k} \epsilon(\lambda )
\exp (i \pi  \trace (S(\lambda )Z )) $.
\end{defi}

By \cite{BFW} the form 
$\BFW(\Lambda, p)$ is a linear combination of Siegel theta-series
of lattices in the genus of $\Lambda $:
For any $k$-dimensional totally  isotropic subspace $E$ of 
$\Lambda / p\Lambda $ let $\Gamma (E) := \langle E , p\Lambda \rangle $
be the full preimage of $E$.
Dividing the scalar product by $p$, one
obtains a lattice $\ ^{1/p} \Gamma (E) := (\Gamma (E)  , \frac{1}{p} (,) ) \in {\cal G}$.
Then we define 
$$\Per (\Lambda ,p):= \sum _{E} \epsilon(E) [ \ ^{1/p} \Gamma (E) ] \in {\cal V}$$
where the sum runs over all $k$-dimensional totally isotropic
subspaces of $\Lambda / p\Lambda $.
As $\epsilon  $ is only defined up to a sign, also $\Per (\Lambda ,p) $
is only well defined up to a factor $\pm 1 $. 
It is shown in \cite[Theorem 2]{BFW} that 
$$\Theta ^{(k)} (\Per (\Lambda ,p)) = \BFW (\Lambda , p ).$$
In analogy to the notation in \cite{KV} we call 
$\Per (\Lambda , p)$ the {\em perestroika} of $\Lambda $.
Clearly $\BFW (\Lambda ,p )$ is in the kernel of the $\Phi $-operator
and hence a cusp form, if the level of $\Lambda $ is
square-free by Theorem \ref{cusp}. 

\subsection{Hecke-actions}\label{Hecke}

Strongly related to the Borcherds-Freitag-Weissauer construction
are the Hecke operators $T(p)$ which define self-adjoint linear
operators on ${\cal V}$ and whose action on theta series
coincides with the one of $T(p)$ in 
\cite[Theorem IV.5.10]{Freitag} and  \cite[Proposition 1.9]{Yoshida} 
up to a scalar factor (depending on the degree of the theta series).
Assume that the genus ${\cal G}$ consists of even $2k$-dimensional
lattices of level $l$.
For primes $p$ not dividing $l$ 
we define $T(p) : {\cal V} \to {\cal V}$ by 
$$T(p) ([\Lambda ]):= \sum _{E}  [ \ ^{1/p} \Gamma (E) ] $$
where the sum runs over all $k$-dimensional totally isotropic
subspaces of $(\Lambda /p\Lambda , Q_p)$.
Note that $T(p)$ is 0 if  
$(\Lambda /p\Lambda , Q_p)$ is not isomorphic to the
sum of $k$ hyperbolic planes.

The following operators commute with the $T(p)$ 
and are usually easier to calculate using the
 Kneser neighbouring-method (see \cite{Kneser}):
For a prime  $p$ 
define  the linear operator $K_p$  by
$$K_p ([ \Lambda ] ) := \sum _{\Gamma } [ \Gamma ], \fa \Lambda \in {\cal G}$$
where the sum runs over all lattices $\Gamma $ in ${\cal G}$
 such that 
the intersection $ \Lambda \cap \Gamma$ has index $p$ in 
$\Lambda $ and in $\Gamma $.
If $p$ does not divide the level  $l$
\cite[Proposition 1.10]{Yoshida} shows that the operators 
$K_p$ are essentially the Hecke operators $T^{(m-1)}(p^2)$ 
(up to a summand, which is a multiple of the identity and a
scalar factor).
Also if $p$ divides $l$, the operators $K_p$ are
self-adjoint:
For  $\Lambda $ and $\Gamma $  in ${\cal G}$, the 
number $n(\Gamma,[\Lambda ])$ 
of neighbours of $\Gamma $ that are isometric to
$\Lambda $ 
equals the number of rational matrices $X\in \GL _{2k}(\Z ) \diag (
p^{-1},1^{2k-1},p) \GL_{2k }(\Z )$ solving 
$$I(\Gamma , \Lambda ): \ \ 
X F_{\Gamma } X^{tr} = F_{\Lambda }$$ (where $F_{\Gamma }$ 
and $F_{\Lambda }$ denote fixed Gram matrices of $\Gamma $ respectively
$\Lambda $) divided by the order of the automorphism group of
$\Lambda $ (since one only counts lattices, $X$ and $gX$ 
have to be identified for all $g\in \GL _{2k}(\Z)$ with 
$g F_{\Lambda } g^{tr} = F_{\Lambda }$).
Mapping $X$ to $X^{-1}$ gives a bijection between the 
set of solutions of $I(\Gamma , \Lambda )$ and 
$I(\Lambda , \Gamma )$. 
Therefore 
$$n(\Gamma , [\Lambda ])  |\Aut (\Lambda )| = 
n(\Lambda , [\Gamma ])  |\Aut (\Gamma )|  .$$

Hence the linear operators $K_p$ and $T(p)$  generate a commutative subalgebra 
$${\cal H}:= \langle T(q), K_p \mid q,p \mbox{ primes }, q \teiltnicht l \rangle \leq \End^s ({\cal V}) $$
of the space of self-adjoint endomorphisms of ${\cal V}$
and ${\cal V}$ has an orthogonal basis 
$(d_1,\ldots , d_h)$,
consisting of common eigenvectors of ${\cal H}$.

For each $1\leq i \leq h$ we define $v(i)\in \{ -1,\ldots, 2k-1 \}$ by
$d_i \in {\cal V}_{v(i)},\ \ d_i \not\in {\cal V}_{v(i)+1} .$
\\
Analogously let $w(i) \in \{ 0,\ldots ,2k \}$ be defined by
$d_i \in {\cal W}_{w(i)},\ \ d_i \not\in {\cal W}_{w(i)-1} .$

\begin{lemma}(\cite[Lemma 2.5]{SieMod})\label{mult1}
Let $1\leq i \leq h$ and assume that $d_i$ generates a full
 eigenspace of ${\cal H}$.
 Then
 $w(i) = v(i)+1$.
 \end{lemma}

If the genus ${\cal G}$ is strongly modular of level $l$, by which we mean that
$\sqrt{d} \Lambda ^{\# , d} \in {\cal G}$ for all $\Lambda \in {\cal G}$ 
and all exact divisors $d$ of $l$, then
the Atkin-Lehner involutions 
$$ W_d: [\Lambda ] \mapsto [\sqrt{d} \Lambda ^{\# ,d } ] $$
for exact divisors $d$ of $l$
define further self-adjoint linear operators on ${\cal V}$.
In this case let 
$$\hat{{\cal H}}:= \langle {\cal H}, W_{d} \mid d \mbox{ exact divisor of } l \rangle . $$

If all lattices in ${\cal G}$ are strongly modular, then $W_d = 1 $ for all 
$d$ and 
 $\hat{\cal H} = {\cal H}$ is commutative.

Again, the Hecke action is rational on ${\cal V}_{\Q }$ hence 
the $\Q$-algebras ${\cal H}_{\Q }$ and $\hat{\cal H}_{\Q}$
spanned by the $K_p$ respectively the $K_p$ and  $W_d$ 
act on ${\cal V}_{\Q }$.

\begin{remark}
For $v\in \{ -1,0,\ldots , 2k-1 \} $ let 
$${\cal D}_v := \langle d_i \mid v(i) = v \rangle .$$
If all eigenspaces of ${\cal H}$ are 1-dimensional, the
 decomposition ${\cal V} = \oplus _{v=-1}^{2k-1} {\cal D}_v$
is preserved under any semi-simple algebra ${\cal A}$
with ${\cal H} \leq {\cal A} \leq \End (V) $ that respects the
filtration.
\end{remark}

\section{Results}\label{results}

The explicit calculations are performed in MAGMA (\cite{MAGMA}).
Fix $l\in {\cal N}$, let  ${\cal G}:= {\cal G}(\Gamma _l)$, ${\cal V} = {\cal V}({\cal G})$  and
denote by $\Lambda _1 := \Gamma _l, \Lambda _2,\ldots, \Lambda _h$ 
representatives of the isometry classes of lattices in ${\cal G}$.
We find that in all cases ${\cal H}   = \langle K_2,K_3 \rangle \cong \C ^h$ 
is a maximal commutative subalgebra of $\End ({\cal V}) $.
Therefore the common eigenspaces are of dimension one
and it is straightforward to calculate an explicit
orthogonal basis
$(d_1,\ldots , d_h)$ of ${\cal V}$ consisting of eigenvectors of ${\cal H}$.
In particular  $v(i)  = w(i) -1 $ for all $i=1,\ldots, h$
by Lemma \ref{mult1}.
Here we  choose 
$d_1 := \sum _{i=1} ^h | Aut(\Lambda _i)|^{-1} [\Lambda _i ]  
 \in {\cal V}_0 \setminus {\cal V}_1$ to be the unit element of ${\cal V}$ 
and  (for $l\neq 7$) $d_h  = \Per (\Gamma _l , p) \in {\cal V}_{k_l-1}$,
where $p$ is
the smallest prime not dividing $l$.
We then determine some Fourier-coefficients of the
series $\Theta ^{(n)} (d_i)$ $(n=0,1,\ldots, k_l)$ to get upper
bounds on $v(i)$.
In all cases the degree-$k_l$ Siegel theta-series of the lattices are
linearly independent hence ${\cal V}_{k_l} = \{ 0 \}$.
Moreover ${\cal V}_{{k_l}-1} = \langle d_h \rangle $ if $l\neq 7$.
We also know that $w(1)=0$ and we may choose $d_2$ such that 
$w(2) =1 $.
By 
 Lemma \ref{mult} and \ref{mult1} 
the product 
$d_j \circ d_i $ lies in $ {\cal W} _{w(i)+w(j)}$.
If the coefficient of $d_h$ in the product is
non-zero, this yields lower bounds on the sum $w(i)+w(j)$
which often  yield sharp lower bound for $w(i)$ and $w(j)$.
The method is illustrated in \cite[Section 3.2]{SieMod} and
an example is given in Paragraph \ref{BW}.

\subsection{The genus of the Barnes-Wall lattice in dimension 16.}\label{BW}

The lattices in this genus are given in \cite{BW}.
The class number is $h=24$ and we find
$$\langle K_2 , K_3 \rangle = {\cal H}_{\Q } \cong \Q ^{13}  \oplus F_1 \oplus F_2 \oplus F_3 $$
where the totally real number fields $F_i \cong \Q[x] / (f_i(x)) $ are given by
$$\begin{array}{ll} 
f_1    = &  x^3 - 11496x^2 + 41722560x - 47249837568 \\
f_2    =  & x^3 - 1704x^2 + 400320x + 173836800  \\
f_3    =  & x^5 - 11544x^4 + 42868800x^3 - 53956108800x^2 + 1813238784000x 
\\ & + 20094119608320000 
\end{array}
$$
and $\langle K_2 , K_3, W_2 \rangle =
\hat{{\cal H}}_{\Q } \cong \Q^{13} \oplus \Mat _3 ( \Q) \oplus \Mat _3(\Q ) \oplus \Mat _5 (\Q ).$
Let $\alpha _i$, $\beta _i$ and $\gamma _j$ ($ i=1,\ldots 3, j=1,\ldots, 5$)
denote the complex roots of 
the polynomials $f_1$, $f_2$ respectively $f_3$.
Let $\epsilon _i$ ($i=1,\ldots , 3$) denote the primitive idempotents of
${\cal H}_{\Q }$ with ${\cal H}_{\Q} \epsilon _i \cong F_i$.

Since the image of ${\cal V}_{\Q }$ under $\Theta ^{(n)}$  has rational
Fourier-coefficients, the functions $v$ and $w$ are constant on the
eigenspaces $E_i = {\cal V} \epsilon _i $ ($i=1,2,3$).
We therefore give their values in one line in the following tabular:

\begin{theorem}
The functions $v$ and and the eigenvalues of $ev_2$ and $ev_3$ of
$K_2$ respectively $K_3$ on $(d_1,\ldots , d_{24})$ 
are as follows: \\
\begin{center}
\begin{tabular}{|l|c|r|r|c|l|c|r|r|}
\hline
$i$ & $v(i)$ & $ev_2$ & $ev_3$ & & 
$i$ & $v(i)$ & $ev_2$ & $ev_3$  \\
\hline
$1$ & $-1$ & $34560$ & $7176640$ & & $15$ & $3,4$ & $1320$ & $8640$ \\
$2$ & $0$ & $ 16200$ & $2389440$ & & $E_2$ & $4$ & $\beta _j$ & $31680$ \\
$3$ & $1$ & $8760$ & $792000$ & & $19$ & $3,4,5$ & $1080$ & $-45120$ \\
$4$ & $1$ & $7128$ & $804288$ & & $20$ & $3,4,5$ & $312$ & $4032$ \\
$E_1$ & $2$ & $\alpha _j$ & $266688$ & & $21$ & $5$ & $-216$ & $8640$  \\
$8$ & $3$ & $2664$ & $90048$ & & $22$ & $5$ & $-216$ & $20928$ \\
$9$ & $3$ & $1320$ & $77760$ & & $23$ & $6$ & $-936$ & $13248$ \\
$E_3$ & $3$ & $\gamma _j$ &$ 100800$ & & $24$ & $7$ & $-2160$ & $39360$ \\
\hline
\end{tabular}
\end{center}
For the dimensions of ${\cal D}_v$ one finds
$$
\begin{array}{|l|ccccccccc|}
\hline
v & -1 & 0 & 1 & 2 & 3 & 4 & 5 & 6 & 7 \\
\hline
\dim ({\cal D}_v) & 1 & 1 & 2 & 3 & 7\mbox{-}10 & 3\mbox{-}5 & 2\mbox{-}4 & 1 & 1 \\
\hline
\end{array}
$$
\end{theorem}

\underline{Proof.}
By explicit calculations of the Fourier-coefficients the values
given in the table are upper bounds for the $v(i)$.
By Lemma \ref{mult1} they also 
provide upper bounds on the $w(i) = v(i) + 1 $.

We see that
$$d_i \circ d_j = A_{ij} d_{24} + \sum _{m=1}^{23} b_{ij}^m d_m $$
with a nonzero coefficient $A_{ij}$ for the following
pairs $(i,j)$:
$$(23,2),\ (22,3),\ (21,4),\ (E_1,E_2),\ (E_3,E_3),\ (8,8), \ (9,9)  $$
(where $(E_1,E_2)$ means that there is some vector in $E_1$ and
some in $E_2$ such that this coefficient is non-zero, similarly $(E_3,E_3)$).
Since $d_m \in {\cal W}_7$ for all $m\leq 23$ and 
$d_j\circ d_i \in {\cal W}_{w(i)+w(j)}$ the inequality 
$w(i)+w(j) \leq 7$ together with $A_{ij} \neq 0$
implies that $d_{24} \in {\cal W}_7$ which is a contradiction.
Hence $w(i) + w(j) \geq 8$ for all pairs $(i,j)$ above.
This yields equality for all values $v(i)$ and $v(j)$ 
for these pairs.
Similarly we get $3\leq v(i) $ for $i=15,19,20$ since $A_{i,i} \neq 0$
for these $i$.
\hfill{q.e.d.}

\begin{Conjecture}
$v(19) = 5$ and  $v(20) = 5 $.
\end{Conjecture}

Since $d_{15}\circ d_2 = \sum _{m=1}^{18} c_m d_m + A_1 d_{19} + A_2 d_{20} $ with $A_1 \neq 0 \neq A_2$, we get 
$w(15) + 1 \geq \max (w(19),w(20))$.

\begin{remark}
If the conjecture is true, then 
$v(15) = 4$ and 
$\dim ({\cal D}_3 )= 7$,
$\dim ({\cal D}_4 )= 4$, and
$\dim ({\cal D}_5 )= 4$.
\end{remark}

Using the formula in \cite[Korollar 3]{Krieg} (resp. \cite[Proposition 1.9]{Yoshida}) we may calculate the eigenvalues of $T^{(m-1)}(3^2)$ from the
one of $K_3$ and compare them with the ones given in \cite[formula (7)]{BK}.
The result suggests that 
 $\Theta ^{(2)}(d_4)$,
$\Theta ^{(4)}(v)$ (for some $v\in E_3$),
 $\Theta ^{(6)}(d_{19})$ and 
 $\Theta ^{(8)}(d_{24})$ are
generalized Duke-Imamoglu-Ikeda-lifts (cf. \cite{Ikeda}) of the elliptic cusp forms
$\delta _8 \theta _{D_4}^i$  (i=3,2,1,0) where
$\delta _8 = \frac{1}{96} (\theta _{D_4}^4-\theta _{\Gamma _3})$ 
is the cusp form of $\Gamma _0(2)$ of
weight 8 and $\theta _{D_4}$ the theta series of the 4-dimensional
2-modular root lattice
$D_4$.  
This would imply that $v(19)=5$ and, with Lemma \ref{mult}, $v(15)=4$.

\subsection{The genus of the Coxeter-Todd lattice in dimension 12.}

For $l=3$ one has $h=10$,
all lattices in this genus are modular, and
${\cal H} _{\Q } = \langle K_2 \rangle \cong \Q ^{10} = \hat{\cal H} _{\Q }$ 

\begin{theorem}
There is some $a\in \{ 0,1 \}$ such that 
the function $v$ and the eigenvalues $ev_2$ of $K_2$ and $e_2$ of $T(2)$ 
are as follows:
\begin{center}
\begin{tabular}{|c|c|r|r|c|c|c|r|r|}
\hline
$i$ & $v(i)$ & $ev_2$ & $e_2$ & & $i$ & $v(i)$ & $ev_2 $ & $e_2$ \\
\hline
$1$ & $-1$ & $2079$ & $151470 $  & & 
$6$ & $3-a$ & $234$ & $7560 $ \\
$2$ & $0$ & $1026$ & $-27540 $  & & 
$7$ & $3$ & $126$ & $2376 $ \\
$3$ & $1$ & $594$  & $17820 $ & &  
$8$ & $3$ & $-36$ & $432 $ \\
$4$ & $1$ & $432$  & $3240 $ & &  
$9$ & $4$ & $-144$ & $-864 $ \\
$5$ & $2$ & $288$  & $-5400 $ & &  
$10$ & $5$ & $-378$ & $1944 $ \\
\hline
\end{tabular}
\end{center}
For the dimensions of ${\cal D}_v$ one finds
$$
\begin{array}{|l|ccccccc|}
\hline
v & -1 & 0 & 1 & 2 & 3 & 4 & 5  \\
\hline
\dim ({\cal D}_v) & 1 & 1 & 2 & 1+a & 3-a & 1  & 1 \\
\hline
\end{array}
$$
\end{theorem}

We conjecture that $a=0$ but cannot prove this using 
Lemma \ref{mult}.

The eigenvalues of $T(2)$ suggest that $\Theta ^{(2)}(d_3)$,
$\Theta ^{(4)}(d_6)$ and $\Theta ^{(6)}(d_{10})$ are
generalized Duke-Imamoglu-Ikeda-lifts (cf. \cite{Ikeda}) of the elliptic cusp forms
$\delta _6 \theta _{A_2}^2$, 
$\delta _6 \theta _{A_2}$,  respectively
$\delta _6 $, where 
$\delta _6 = \frac{1}{36} (\theta _{A_2}^6-\theta _{\Gamma _3})$ 
is the cusp form of $\Gamma _0(3)$ of
weight 6 and $\theta _{A_2}$ the theta series of the hexagonal lattice
$A_2$.  
This would imply $v(3) = 1$, $v(6)=3$ and $v(10)=5$ and hence $a=0$.

\subsection{The genus of the 5-modular lattices in dimension 8.}

The class number of this genus is $h=5$,
all lattices in this genus are modular, and
${\cal H} _{\Q } = \langle K_2 \rangle \cong \Q ^5 = \hat{\cal H} _{\Q }$ 

\begin{theorem}
For $l=5$ one has 
$\dim ({\cal D}_{v} ) =  1 $ for $v=-1,0,1,2,3$. 
The function $v$ and the eigenvalues $ev_2$ of $K_2$ and 
$e_p$ of $T(p)$ ($p=2,3$) 
are  given in the following table:
\begin{center}
\begin{tabular}{|c|c|r|r|r|c|c|c|r|r|r|}
\hline
$i$ & $v(i)$ & $ev_2$ & $e_2$ & $ e_3$ & & $i$ & $v(i)$ & $ev_2 $ & $e_2$ & $ e_3$ \\
\hline
$1$ & $-1$ & $135$ & $270$ & $2240 $ & & 
$4$ & $2$ & $-8$ & $-16 $ & $56 $ \\
$2$ & $0$ & $70$  & $-120 $ & $ 160 $ & & 
$5$ & $3$ & $-60$ & $10$ & $ 420 $ \\
$3$ & $1$ & $42$  & $84 $ & $ 256 $ & &  & &  & & \\
\hline
\end{tabular}
\end{center}
\end{theorem}


\subsection{The genus of the strongly 6-modular lattices in dimension 8.}

The class number of  ${\cal G}(\Gamma _6)$ is
$h=8$,  the Hecke-algebras are
$\hat{{\cal H}} _{\Q } = \langle K_2 , W_2 \rangle \cong \Q^5 \oplus \Mat _3(\Q ) $
and  ${\cal H}_{\Q }= \langle K_2 \rangle  \cong \Q^5 + \Q[x]/(f(x))$
where $$f(x) = x^3 - 66x^2 - 216x + 31104.$$ 
Let $\delta _i \in \R $ ($i=1,2,3$) denote the roots of $f$.
\begin{theorem}
Then the function $v$ 
and the eigenvalues $ev_2$ of $K_2$ and $e_5$ of $T(5)$
are given in the following table:
\begin{center}
\begin{tabular}{|c|c|r|r|c|c|c|r|r|}
\hline
$i$ & $v(i)$ & $ev_2$  & $ e_5 $ & & $i$ & $v(i)$ & $ev_2$ & $ e_5 $  \\
\hline
$1$ & $-1$ & $144$   & $39312 $ & & $E$ & $1$ & $\delta _j $  & $ 1872 $ \\
$2$ & $0$ & $54$  & $1872 $ & 
& $7$ & $2$ & $-6$  & $432 $ \\
$3$  & $1$ & $ 18 $   & $ 1008 $ & & 
$8$ & $3$ & $-36$   & $ 4752 $ \\
\hline
\end{tabular}
\end{center}
Hence 
$\dim ({\cal D}_v) = 1$ for $v=-1,0,2,3$ and $\dim({\cal D}_1) = 4$.
\end{theorem}

\subsection{The genus of the 7-modular lattices in dimension 6.}

The class number is $h=3$, all lattices are modular, and 
$\hat{\cal H} _{\Q } = {\cal H} _{\Q } = \langle K_2 \rangle \cong \Q^3$.
In contrast to the other genera,
the perestroika $\Per (\Gamma _7,2) $ and hence also $\BFW (\Gamma _7,2)$ 
vanishes due to the fact that the image of $\Aut (\Gamma _7)$
in $GO_6^+(2)$ is not contained in the derived subgroup $O_6^+(2)$.
In fact, $\Theta ^{(2)} $ is already injective. 
Since the discriminant of the space is not a square modulo $3$ and $5$,
the Hecke operators $T(3)$ and $T(5)$ vanish.

\begin{theorem}
We have $v(i) = i-2$ for $i=1,2,3$ and hence 
$\dim ({\cal D}_v) = 1 $ for $v=-1,0,1 $.
The eigenvalues of $K_2$ are $35$, $19$, and $5$,
the ones of $T(2)$ are $30$, $-18$, and $10$, and
$T(11)$ has eigenvalues 
$2928$, $ -144$, and $ 248$.
\end{theorem}

\subsection{The genus of the strongly $l$-modular lattices in dimension 4
for $l=11,14,15$.}

For $l=11,14,15$ the genus ${\cal G}(\Gamma _l)$ consists of 3 isometry classes
and ${\cal H} _{\Q } = \langle K_2 \rangle \cong \Q ^3 = \hat{\cal H} _{\Q }$ since all lattices in the genus
are strongly modular.
\begin{theorem}
For $l=11,14,15$ one has 
$\dim ({\cal D}_{v}) =  1 $ for $v=-1,0,1 $. 
The eigenvalues $ev_2$ of $K_2$ and $e_p$ of $T(p)$ for primes $p\leq 7$ not dividing
$l$
are  given in the following table:
\begin{center}
\begin{tabular}{|c|c||r|r|r|r|r||r|r|r||r|r|r|}
\hline
& & \multicolumn{5}{|c||}{$l=11$} &
\multicolumn{3}{|c||}{$l=14$}  &
\multicolumn{3}{|c|}{$l=15$}  \\
\hline
$i$ & $v(i)$ & $ev_2$ & $e_2$ & $e_3$ & $e_5$ & $e_7$ &
 $ev_2$ & $e_3$ & $e_5$  &
 $ev_2$ & $e_2$ & $e_7$  \\
\hline
$1$ & $-1$ & $9$ & $6$ & $8$ & $12$ & $16$  &
 $8$ & $8$ & $12$ &
 $9$ & $6$ & $16$  \\

$2$ & $0$ & $4$ & $-4$ & $-2 $ & $2 $ & $-4 $  &
 $2$ & $-4$ & $0$  &
 $1$ & $ -2 $ & $ 0 $ \\

$3$  & $1$ & $ -6 $ & $1$ & $3 $ & $ 7 $ & $ 6 $ &
 $ -4 $ & $2$ & $6$ &
 $ -3 $ & $2$ & $8 $ \\
\hline
\end{tabular}
\end{center}
\end{theorem}

\subsection{The genus of the $23$-modular lattices in dimension 2.}
In the smallest possible dimension $2$ the genus ${\cal G}(\Gamma _{23})$ consists of only
2 isometry classes and ${\cal H} _{\Q } = \langle K_2 \rangle \cong \Q ^2 =
\hat{\cal H} _{\Q }$ for the same argument that all lattices in the genus are
modular.

\begin{theorem}
For $l=23$ one has
$\dim ({\cal D}_{v}) =  1 $ for $v=-1,0$.
One has $v(1) = -1$, $v(2)= 0$, 
$d_1 K_2  = 2 d_1 $ and
$d_2 K_2  = - d_2 $.
For the $T(p)$ for primes $p<23$ we find
$T(2)=T(3)=T(13) = K_2$ and $T(5)=T(7)=T(11)=T(17)=T(19)=0$.
\end{theorem}

\noindent Authors' addresses:\\
Gabriele Nebe,
Lehrstuhl D f\"ur Mathematik, RWTH Aachen,
 52056 Aachen, Germany, nebe@math.rwth-aachen.de
\\
Maria Teider,
Abteilung Reine Mathematik, Universit\"at Ulm, 89069 Ulm, Germany
maria.teider@mathematik.uni-ulm.de
\end{document}